\newif\ifsmfart
 \IfFileExists{smfart.cls}
   {\documentclass[12pt,english]{smfart}\smfarttrue}
   {\message{^^J*** This file should be typeset with smfart.cls ***^^J^^J}
    \documentclass[12pt]{amsart}\let\Subsection\subsection}

\usepackage{amssymb}

\IfFileExists{mathrsfs.sty}
  {\usepackage{mathrsfs}\let\mathcal\mathscr}
  {\let\mathscr\mathcal}

\advance\headheight 2pt
\calclayout

\numberwithin{equation}{section}
\makeatletter

\let\c@equation\c@subsection
\def\theequation{\thesection.\arabic{equation}}

\def\l@table{\@tocline{0}{3pt plus2pt}{0pt}{}{\itshape}}

\makeatother

\theoremstyle{plain}
\newtheorem{prop}[subsection]{Proposition}
\newtheorem{thm}[subsection]{Theorem}
\newtheorem{cor}[subsection]{Corollary}

\newtheorem{lem}[subsection]{Lemma}

\theoremstyle{definition}

\theoremstyle{remark}
\newtheorem{rem}[subsection]{Remark}
\newtheorem{rems}[subsection]{Remarks}
\newtheorem{exam}[subsection]{Example}

\def\A{{\mathbf A}}
\def\C{{\mathbf C}}
\def\F{{\mathbf F}}
\def\G{{\mathbf G}}
\def\N{{\mathbf N}}
\def\P{{\mathbf P}}
\def\Q{{\mathbf Q}}
\def\R{{\mathbf R}}

\def\Z{{\mathbf Z}}

\let\ra\rightarrow

\let\epsilon\varepsilon \let\eps\epsilon
\let\epsilon\varepsilon
\let\phi\varphi

\let\leq\leqslant
\let\geq\geqslant
\let\le\leq
\let\ge\geq

\def\eff{\text{\upshape eff}}
\def\abs#1{\left\lvert{#1}\right\rvert}
\def\norm#1{\left\|{#1}\right\|}

\def\DeclareMathOperator#1#2{\def #1{\operatorname{#2}}}
\DeclareMathOperator{\Re}{Re} 
\DeclareMathOperator{\Pic}{Pic}
\DeclareMathOperator{\Br}{Br}
\DeclareMathOperator{\Spec}{Spec}
\DeclareMathOperator{\Proj}{Proj}
\DeclareMathOperator{\rk}{rk}
\DeclareMathOperator{\codim}{codim}
\DeclareMathOperator{\div}{div}
\DeclareMathOperator{\Val}{Val}
\DeclareMathOperator{\vol}{vol}

\def\dx{\,{\mathrm d\mathbf x}}

\title[Points of bounded height]
      {Points of bounded height \\
       on equivariant compactifications \\ of vector groups, I}

\author{Antoine Chambert-Loir}
\address{
Institut de math\'ematiques de Jussieu, Boite 247 \\ 
4, place Jussieu \\ 
F-75252 Paris Cedex 05 }
\email{chambert@math.jussieu.fr}

\author{Yuri Tschinkel}
\address{Dept. of Mathematics, U.I.C.\\
Chicago, (IL) 60607-7045,  U.S.A. }
\email{yuri@math.uic.edu}

\begin{document}
\def\smfandname{\&}
\date{\today}

\begin{abstract}
We prove asymptotic formulas  for the number
of rational points of bounded height
on certain equivariant compactifications of the affine plane.
\end{abstract}
 
\ifsmfart 
 \begin{altabstract}
 Nous \'etablissons un d\'eveloppement asymptotique
 du nombre de points rationnels de hauteur born\'ee sur
 certaines compactifications \'equivariantes du plan affine.
 \end{altabstract}
\fi

\maketitle

\tableofcontents

\section*{Introduction}
\def\thesubsection{\arabic{subsection}}

Let $X$ be a smooth projective algebraic variety 
defined over a number field $F$ and $X(F)$ the 
set of rational points of $X$. 
Let ${\mathcal L}$  be a metrized ample line bundle 
and 
\[
H_{\mathcal L} : X(F)\ra \R_{>0}
\]
the associated exponential 
height (cf.  \cite{peyre99}, \cite{salberger99}). 
We are interested
in the asymptotic behaviour of the counting function
\[
N(U,{\mathcal L},B):=
\# \{ x\in U(F) \,|\, H_{\mathcal L}(x) \le B \}
\]
as $B\ra \infty$, where $U\subset X$ is some Zariski open 
subset. 
Batyrev, Manin~\cite{batyrev-m90} and
Peyre~\cite{peyre95} described
 a conjectural framework for such asymptotics
for varieties with ample (or, more generally, effective) 
anticanonical class
(see also~\cite{franke-m-t89}, \cite{peyre99}
and references therein).
In particular, it is expected that for appropriate $X$ and Zariski
open subsets $U\subset X$,
\[
N(U,K_X^{-1},B)= 
\frac{\Theta(X)}{r!}B(\log B)^{r}(1+o(1)),
\]
as $B\ra \infty$.
Here  $K_X^{-1}$ is the metrized anticanonical 
line bundle on $X$,  $r=\rk \Pic(X)-1$ and $\Theta(X)$
is a product of a Tamagawa type number $\tau(K_X)$
(which depends on the metrization), a rational number $\alpha (X)$
defined in terms of the cone of effective 
divisors $\Lambda _{\eff}(X)$    
and the order $\abs{\Br(X)/\Br(F)}$ of the non-trivial part of the
Brauer group. There is a similar description for 
arbitrary polarizations ${\mathcal L}$ (Batyrev, Tschinkel,
cf. \cite{batyrev-t99}). 

These conjectures have been proved for
flag varieties (\cite{franke-m-t89}),
toric varieties (\cite{batyrev-t99}),
and toric bundles induced from torsors  
(\cite{strauch-t99}).
The proofs use a precise {\em combinatorial} description 
of all geometric and arithmetic invariants of 
the varieties: line bundles, 
metrizations of the line bundles etc. 
(for example, in terms of lattices, cones and fans). 
Such a description is possible because
representations of reductive groups are rigid (don't admit
deformations). Consequently, the 
corresponding varieties don't have moduli.

The only other known approach to asymptotics of 
rational points on algebraic varieties is 
the classical circle method in analytic number theory.
Varieties which can be treated by this method
do admit moduli. However, one of the drawbacks
is that so far it works only for varieties
which are  complete intersections of small degree
$d$ in projective spaces $\P^n$ of large dimension
(very roughly, $n\gg 2^d$).
In particular, these complete intersections
have Picard group~$\Z$.  
There is a promising generalization of
the circle method to complete intersections in 
other varieties (for example, toric varieties) due
to E. Peyre (cf. \cite{peyre98}), which 
should provide examples of asymptotics for varieties
with moduli and with $\Pic(X)$ of higher ranks, once
the necessary estimates are established.  
As a reference to the circle method 
let us mention the papers by 
H.~Davenport, R.~Heath-Brown, C.~Hooley on smooth cubic 
hypersurfaces (cf. \cite{davenport63}, \cite{heath-brown83}, 
\cite{hooley88}), by B.~Birch and by W.~Schmidt 
on general complete intersections (cf. \cite{birch62}, 
\cite{schmidt85}).

In this paper we prove
asymptotics of rational points of bounded height
on varieties which admit moduli and which 
at the same time are closely related to (non-reductive)
linear algebraic groups.
More precisely, we consider equivariant compactifications
of the additive group $\G_a^n$.
For $n=2$ it can be shown that all such 
compactifications are obtained as blow-ups of $\P^2$, 
or Hirzebruch surfaces in points 
which are fixed under the action of $\G_a^2$. 
Notice that a variety (even $\P^2$) may 
admit non-isomorphic structures as an equivariant 
compactification of $\G_a^n$.  
A similar ``minimal model program'' 
of equivariant compactifications of ${\mathbf G}_a^n$ 
is a non-trivial problem already for $n=3$ 
(see~\cite{hassett-t98}).

In this paper we study in detail the example of  
a blow-up of $\P^2$ in $r$  $\Q$-rational points which 
are all contained in the line at infinity
$\P^1$ (with the equation $x_0=0$). 
The moduli space of
such surfaces $X$ is ${\mathcal M}_{0,r}$. 
It is easy to see that $X$ is a smooth projective
equivariant compactification of $\G_a^2$
with $\Pic(X)=\Z^{r+1}$, trivial Brauer group
and a simplicial cone of effective divisors $\Lambda _{\eff}(X)$.  
Denote by $U\simeq \G_a^2\subset \P^2$ the 
complement to $x_0=0$. 
Then for $\Re(s)\gg 0$, the series 
\[
Z(U,K_X^{-1},s)=\sum_{{\mathbf x}\in U(\Q)}
H_{K_X^{-1}}({\mathbf x})^{-s}
\] 
is absolutely and uniformly convergent and defines a holomorphic
function.
One of the main results of this paper is the following:

\begin{thm}\label{theo:main}
There exists a function $h$ which is holomorphic in
the domain $\Re(s)>2/3$ such  that 
\[
Z(U,K_X^{-1},s)=\frac{h(s)}{(s-1)^{r+1}}
 \quad\text{and}\quad
h(1)=\alpha (X)\tau(K_X)\neq 0.  \]
\end{thm}

A standard Tauberian theorem implies that $X$ satisfies
Peyre's refinement of Manin's conjecture:
\begin{cor}
We have the following asymptotic formula:
\[
N(U,K_X^{-1}, B)= \frac{\alpha (X)\tau(K_X)}{r!} 
B(\log B)^{r}(1+ o(1)),
\]
as $B\ra \infty$.
\end{cor}

In fact, we will prove asymptotics for every ${\mathcal L}$ on $X$
such that its class is contained in the interior of $\Lambda _{\eff}(X)$. 
We will also give estimates for the growth
of the function $h(s)$ in vertical strips in 
the neighborhood of $\Re(s)=1$. 
It is well known that this implies a more precise asymptotic expansion
for the counting function $N(U,{\mathcal L}, B)$,
see Corollary~\ref{coro:last}.

\bigskip

We will address the compactifications of $\G_a^n$ (with $n>2$)
in subsequent papers.

{\itshape Acknowledgements.} 
The work of the second author was partially supported by
the NSA.

\def\thesubsection{\thesection.\arabic{subsection}}

\section{Geometry}
\label{sect:geometry}

\subsection{Generalities}
Let $G$ be an algebraic group and $X$ a smooth projective 
variety with an action of $G$.
We denote by $\Pic^G(X)$ the (abelian) group of isomorphy classes of
$G$-linearized line bundles on $X$ (cf.~\cite{mumford-f-k94}, Chap.~1,
\S\,3, Def.~1.6).
We shall say that the variety $X$ is an equivariant compactification of
$G$ if $X$ has an open dense subset $U$ which is equivariantly isomorphic
to $G$.
Well known examples are given by toric varieties 
which are equivariant compactifications of tori (algebraic
groups isomorphic to $\G_m^n$ over the algebraic closure of $F$,
where $\G_m=\Spec(F[x,x^{-1}])$ is the multiplicative group scheme).

In this paper, we are interested in equivariant compactifications of
$\G_a^n$, where 
$\G_a=\Spec(F[x])$ is the additive group scheme (we may call
them {\em addic} varieties\footnote{Scherzhafter Vorschlag von 
Yu. I. Manin, Weihnachten 1998}). 

Notice that  a variety can be an equivariant compactification of
a group $G$ in many non-isomorphic ways, as the following example shows.

\begin{exam} 
The projective plane $\P^2$ is an equivariant compactification of
$\G_a^2$ in (essentially) two non-isomorphic ways. One of the possible
actions is the standard translation action, fixing a line 
$\P^1$ at infinity. All 1-parameter subgroups are lines. 
The other action has exactly one fixed point and 
generic 1-parameter subgroups are conics (cf. \cite{hassett-t98}
for more details, esp.~Prop.~3.2). 
\end{exam}

We quote from~\cite{hassett-t98} the following general geometrical
facts about equivariant compactifications of additive groups.
\begin{prop}\label{prop:geometry}
Let $X$ be a smooth projective
equivariant compactification of $\G_a^n$ and 
$D=X\backslash \G_a^n$ the boundary. 
\begin{enumerate}
\item 
The boundary $D$ is a pure codimension 1 subvariety.
\item  
The Picard group $\Pic(X)$ is freely generated
by the irreducible components $D_0,...,D_{r}$ of $D$.
\item
The closed cone of effective divisors $\Lambda _{\eff}(X)\subset
\Pic(X)_{\R}$ is given by
\[
\Lambda _{\eff}(X)= \bigoplus_{k=0}^{r} \R_+ [D_k].
\]
\item There exist integers $n_k>1$ such that
the anticanonical class is given by 
\[
[K_X^{-1}]=\sum_{k=0}^{r} n_k[D_k].
\]
\end{enumerate}
\end{prop}

\subsection{Blow-ups}
The basic example of an equivariant compactification 
of $\G_a^n$ is the projective space $\P^n$, with
$\G_a^n$ acting on $ \P^n=\Proj(F[x_0,\dots,x_n])$
by translation:
\[ ((t_1,\dots,t_n), (x_0:\dots:x_n))\mapsto
         (x_0:x_1+t_1x_0:\dots:x_n+t_nx_0)\]
which stabilizes the ``hyperplane at infinity''
given by the equation $x_0=0$.
In this paper 
we consider blow-ups of the projective space 
$\pi \colon  X\ra \P^n$
in a subscheme $Z\subset \P^n$
of pure codimension $\ge 2$, which is 
contained in this hyperplane.
We denote by $\mathcal I_Z$
the ideal sheaf of $Z$ in $\P^n$ so that
the blow-up is defined by 
the formula $X=\mathbf{Proj}(\bigoplus_j {\mathcal I}_Z^j)$.
As $Z\subset Z_0$ is fixed by $\G_a^n$,
the universal property of the blow-up implies that
the action of $\G_a^n$ on $\P^n$ lifts uniquely to an action
on $X$.

The geometry of blow-ups of arbitrary subschemes
can be very complicated.  We shall assume that $\mathscr I_Z$
is the product of $r$ ideals $\mathscr I_{Z_k}$
($1\leq k\leq r$),
where the $Z_k$ are integral subschemes of the hyperplane
at infinity in $\P^n$.
The universal property of blow-ups says that $X$ is
the universal scheme mapping to $\P^n$ on which $\mathscr I_Z$
becomes invertible. An easy lemma in commutative algebra
implies that on $X$, the $\mathscr I_{Z_k}$ are themselves
invertible. (Blowing up a product of ideals is the
universal way to make these ideals invertible; it is the
same as blowing up successively $Z_1$, then the strict transform
of $Z_2$, etc.) In particular, $\pi$ factors as
$X\ra X_k\ra \P^n$, $X_k$ being the blow-up of $Z_k$ in $\P^n$.

On $X$, we now have Cartier divisors $D_k$ (which are the inverse
images of the $Z_k$) and line bundles $\mathscr O_X(D_k)$
equipped with a canonical section
$\mathsf s_{D_k} \in\Gamma(X,\mathscr O_X(D_k))$.
Moreover, $\mathsf s_{D_k}$ and $\mathscr O_X(D_k)$ are pullbacks
of similar objects on $X_k$ (which we will denote by the same
letters).
Note also that 
$\mathscr O_X(-D_k)=\mathscr I_{Z_k}\cdot\mathscr O_X\subset\mathscr O_X$
and that by definition, $\mathsf s_{D_k}$ is obtained by
dualizing the pull-back of the canonical inclusion 
$\mathscr I_{Z_k}\ra\mathscr O_{\P^n}$.

Let $D_0$ be the strict transform of $Z_0$ in $X$.
We have a canonical isomorphism:
\begin{equation}\label{eqn:z0}
\pi^*{\mathcal O}_X(Z_0)\simeq {\mathcal O}_X(D_0+\sum_{k=1}^r D_k).
\end{equation}
Denoting by $\mathsf s_{D_0}$ the canonical section of $\mathscr O_X(D_0)$,
the tensor product $\bigotimes_{k=0}^r \mathsf s_{D_k}$
equals the pull back on $X$ of the canonical section
of $\mathscr O(Z_0)$.

The $D_k$'s ($0\leq k\leq r$) form a basis of the Picard group of $X$.
We identify the anticanonical sheaf in these coordinates:
\begin{prop} \label{eqn:kx}
If $Z$ is smooth, then
$X$ is a smooth variety and its anticanonical sheaf is given by
\begin{align*}
K_{X}^{-1} &=\pi^*((n+1)Z_0)\otimes \bigotimes_{k=1}^{r}
{\mathcal O}(-(\codim(Z_k)-1)D_k) \\
&=\bigotimes_{k=0}^{r}
{\mathcal O}((2+\dim(Z_k))D_k).
\end{align*}
\end{prop}
\begin{proof}
See Hartshorne~\cite{hartshorne77}, Ex.~8.5, p.~188.
\end{proof}

\subsection{Metrizations on blow-ups}

Let $S$ be the spectrum of a Dedekind ring
(which will be the ring of integers in $F$, or a localization
of it, or a completion) or the spectrum  of a field which is equipped
with a valuation.
Let $X$ be a projective scheme over $S$.
For a locally free sheaf $\mathscr E$
of finite rank on $X$,
there are several notions of metrizations corresponding to these
various cases. We recall briefly the definitions.
\begin{itemize}
\item If $S=\Spec(F)$, the spectrum of a field
endowed with a valuation,
a metric on $\mathscr E$ is a family of norms on 
the fibres $\mathscr E_x$ for $x\in X(\bar F)$, which vary
continuously with $x$.
\item If $\bar F=\C$, one may ask that
the dependence is $\mathscr C^\infty$, and---independently---
that the metrics are hermitian in the fibers.
\item If $S=\Spec(F)$, where $F$ is the fraction
field of a discrete valuation ring $R$,
any flat and projective model $(\mathscr X,\mathscr E)$ over
$\Spec(R)$ determines a metric according to which a section
is of norm $\leq 1$ at a point iff it is integral.
\item If $S=\Spec(F)$, $F$ being a number field,
an adelic metric on $\mathscr E$ is a collection of
metrics for all $X_v/F_v$, $v$ being the different places of $F$.
Moreover, one assumes that there exists a model over 
$\Spec(\mathfrak o_F)$ which gives the same metrics except
at a finite number of places. At these exceptional
places the ratios of the two metrics are assumed to be bounded functions
on $X$.
\end{itemize}
The usual definitions of metrics on subsheafs, quotients, hom's,
etc. are compatible with these notions.

Let $X$ be a quasi-projective flat scheme over $S$,
$\mathscr I$ a sheaf of ideals on $X$ and $Z=V(\mathscr I)$.
Let $\pi:Y\ra X$ be the blow-up of $V(\mathscr I)$,
$Y=\mathbf{Proj}(\bigoplus_n \mathscr I^n)$.
On $Y$, 
the inverse image of $Z$ becomes a Cartier divisor $D$
and the line bundle $\mathscr O(D)$ is
equipped with  a canonical section $\mathsf s_{D}$.
We want to endow $\mathscr O_Y(D)$ with a metric and
to give a formula for the norm of $\mathsf s_D$ at
any point of $Y \backslash |D|$.
Note that $\mathscr O(-D)=\mathscr I\cdot\mathscr O_Y\subset\mathscr O_Y$
and that $\mathsf s_{D}$ is the pull-back of the canonical inclusion 
$\mathscr I\ra\mathscr O_X$.

Choose a locally free sheaf $\mathscr E$ of finite rank on $X$ 
with a section $\sigma_Z\in\Gamma(X,\mathscr E)$ whose scheme of zeroes is $Z$
(existence follows from the quasi-projectivity of $X$).
This induces a surjective homomorphism $\phi:\mathscr E^\vee\ra\mathscr I$
and a closed immersion
$Y\hookrightarrow \P(\mathscr E^\vee)$ such that
$\mathscr O(-D)=\mathscr I\cdot\mathscr O_Y =\mathscr O_{\P}(1)$ and
the universal quotient map
$\pi^*\mathscr E^\vee\ra\mathscr O_{\P}(1)$ on $Y$
is the pullback of $\phi$.
Hence, to metrize $\mathscr O_Y(D)$
it is sufficient to endow $\mathscr E^\vee$ with a metric.
The quotient metric  on $\mathscr O_{\P}(1)$  is defined as follows:
for any local section $s$ of $\mathscr O_{\P}(1)$, we pose
\[ \norm{s}(y) = \inf_{t} \norm{t}(y)  \]
where the infimum is on the local sections $t$ of $\pi^*\mathscr E^\vee$
mapping to $s$ under the canonical surjection
$\pi^*\phi:\pi^*\mathscr E^\vee\ra\mathscr O_{\P}(1)$.

Restrict this to $Y$. This gives a norm on $\mathscr O_{\P}(1)|_{Y}
= \mathscr O_Y(-D)$. The dual norm on  $\mathscr O_Y(D)$ 
of the canonical section $\mathsf s_D$
is given  by the formula
\[ \norm{\mathsf s_D} = \sup_{s\neq 0} \frac{\abs{\langle \mathsf s_D,s\rangle}}
    {\norm{s}}
    =  \sup_t \frac{\abs{\langle \mathsf s_D,\pi^*\phi(t)\rangle}}
            {\norm{t}} \]
the last supremum being over the non-zero local
sections $t$ of $\pi^*\mathscr E^\vee$.
But, away from $D$ on the blow-up, $t$ comes from a local section
of $\mathscr E^\vee$ and
$\langle \mathsf s_D,\pi^*\phi(t)\rangle$
is exactly the image of $t$
under the surjection $\phi:\mathscr E^\vee\ra\mathscr I$.
Hence,
 $\norm{\mathsf s_D}$ is equal to the norm of $\phi$, viewed
as a homomorphism $\mathscr E^\vee\ra\mathscr O_X$,
which by the definition of the dual norm on $\mathscr E^\vee$
is exactly the norm of the original section $\sigma_Z\in\Gamma(X,\mathscr E)$.

(This can be simplified if one regards the blow-up as the closure
of the graph of the map $X\setminus Z\ra \P(\mathscr E^\vee)$
induced by $\sigma_Z$.)

Note the precise meaning of these calculations:
\begin{itemize}
\item they are valid if $S$ is any field with a valuation;
\item if $S$ is a discrete valuation ring, arithmetic intersection
on the integral model gives a result which is compatible with
the metrized theory on the generic fibre if the metric on $\mathscr E$
comes from the model;
\item if $S$ is the ring of integers of a number field, they show
that we have an adelic metric in the sense of Arakelov geometry
provided $\mathscr E$ is equipped with an adelic metric.
\end{itemize}

Hence, we have the following theorem:
\begin{thm}
Let $X$ be an algebraic variety over a field $F$,
$\mathscr I\subset\mathscr O_X$ a sheaf of ideals on $X$
and $\pi:Y\ra X$ the blow-up of $\mathscr Y$.
Let $\mathscr E$ be a locally free sheaf of finite rank on
$X$ with a section $\sigma_Z\in\Gamma(X,\mathscr E)$ such
that $V(\mathscr I)=\div(\sigma_Z)$ as schemes.

Assume $\mathscr E$ is given a metric. Then 
the line sheaf $\mathscr O_Y(D)$
corresponding to the exceptional divisor $D$ on $Y$
has a canonical metric 
such that the norm of its canonical section $\mathsf s_D$
is given by the formula:
\begin{equation}
 \norm{\mathsf s_D}(y) = \norm{\sigma_Z}(\pi(y)). 
\end{equation}
\end{thm}

In particular, if $\mathscr L_1,\dots,\mathscr L_r$ are line bundles on $X$
with sections $\mathsf s_i$ such that, as a scheme,
$Z=\bigcap \div(\mathsf s_i)$,
we may take $\mathscr E=\bigoplus \mathscr L_i$,
$\sigma_Z=(\mathsf s_i)$. Assume the $\mathscr L_i$ to be metrized and endow
$\mathscr E$ with the associated hermitian metric
(resp.\ with the $\ell^\infty$-metric at non-archimedian places).
The preceeding theorem implies that $\mathscr O_Y(D)$ may be
metrized in such a way that
\begin{equation}
 \norm{\mathsf s_D}^2(y) = \sum_{i=1}^r \norm{\mathsf s_i}(\pi(y)). 
\end{equation}

In particular, if $X=\P^n$, $\mathscr L_i= \mathscr O_{\P^n}(n_i)$,
$\mathsf s_i$ corresponds to a homogeneous polynomial $g_i$ of degree $n_i$
and, if $\pi(y)=(x_0:\dots:x_n)$,
\begin{equation}
 \norm{\mathsf s_D}^2 (y)
= \sum_{i=1}^r \frac{ \abs{g_i(x_0,\dots,x_n)}^2}
                    { \big( \sum_{j=0}^n \abs{x_i}^2 \big)^{n_i}}. 
\end{equation}

As a last example, assume that $X=\P^n$ and $Z$ is an integral divisor
in $Z_0$. Then, the homogeneous ideal of $Z$ is of the form
$(x_0,f(x_1,\dots,x_n))$ for some homogeneous polynomial $f$
of degree $d\geq 1$. If $\pi(y)=(1:x_1:\dots:x_n)$, then 
\begin{equation}
 \norm{\mathsf s_D}^2(y)
          = \frac {1}{1+ \sum_{j=1}^n \abs{x_j}^2 }
             + \frac{\abs{ f(x_1,\dots,x_n)}^2}
                    {\big( 1+ \sum_{j=1}^n \abs{x_j}^2 \big)^d}. 
\end{equation}

All these formulas have analogues at non-archimedian places
with the sum of the squares being replaced by their maximum.

\subsection{R\'esum\'e}
Let $F$ be a number field.
For $1\leq k\leq r$, choose a finite family of
homogeneous polynomials $g_{k,j}\in F[x_0,\dots,x_n]$
of degree $d_{k,j}$ generating a prime ideal $\mathcal I_{Z_k}$
corresponding to an integral subscheme $Z_k\subset\P^n$.
Let $\pi:X\ra\P^n$ be the blow-up of the ideal
$\mathscr I =\mathscr I_{Z_1}\cdots \mathscr I_{Z_k}$.
On $X$, the inverse image of $Z_k$ is a Cartier divisor $D_k$
whose associated line bundle $\mathscr O_{X}(D_k)$ can be
adelically metrized
so that the norm of its canonical section $\mathsf s_{D_k}$ at a
point $x\in X$ mapping to $(x_0:\dots:x_n)\in\P^n$
is given by
\[ \norm{\mathsf s_{D_k}}_v(x)
= \max_j \frac{\abs{g_{k,j}(x_0,\dots,x_n)}_v}
              {\max(\abs{x_0}_v,\dots,\abs{x_n}_v)^{d_{k,j}}} \]
at finite places $v$, and by
\[ \norm{\mathsf s_{D_k}}^2_v(x)
= \sum_j \frac{\abs{g_{k,j}(x_0,\dots,x_n)}^2_v}
              {(\abs{x_0}^2_v+\dots+\abs{x_n}^2_v)^{d_{k,j}}} \]
if $v$ is an archimedian place.

We shall henceforth assume that $Z_k$ is contained in the
hyperplane at infinity $x_0=0$. Then one may assume
that one of the $g_{k,j}=x_0$ and that the others do not
depend on $x_0$.
The universal property of the blow-up implies that
$\pi:X\ra\P^n$ is an isomorphism over $\G_a^n\simeq \{x_0\neq 0\}$
and that the action of $\G_a^n$ on $\P^n$ lifts 
to an action on $X$ and to an action on the line
bundles $\mathscr O_X(D_k)$.

The following proposition
can be deduced, either through explicit computations with the formulas
defining $\norm{\mathsf s_{D_k}}$, or by an abstract argument
involving schemes over $\Spec\mathfrak o_F$.

\begin{prop} \label{prop:invariance}
Assume that
the polynomials $g_{k,j}$ have coefficients in $\mathfrak o_F$
and that they 
generate the homogeneous ideal $\mathscr I_{Z_k}\cap\mathfrak
o_F[x_0,\dots,x_n]$.\footnote{This means
that the subscheme $V((g_{k,j})_j)$ of $\P^n_{\mathfrak o_F}$
is projective and flat over $\mathfrak o_F$.}
Then, for each place $v$ of $F$, the standard maximal compact
subgroup of $\G_a^n(F_v)$ acts isometrically on $\mathscr O(D_k)$.
\end{prop}

Let $D_0$ be the strict transform of the hyperplane
at infinity under $\pi$.
The line bundle $\mathscr O_X(D_0)$ is the pull-back
on $X$ of the $\mathscr O_{\P^n}(1)$
and we shall equip it with its standard metric
(given by the formulas above, the family
of $g_{0,j}$ being reduced to $x_0$).
By means of equations
(\ref{eqn:z0}, \ref{eqn:kx}),
we then can metrize the line bundles ${\mathcal O}_X(D_0)$ and 
$K_X^{-1}$.

\section{Heights, Poisson formula}
\label{sect:zeta}

\subsection{Product formula and heights}
We  recall some conventions concerning absolute values
in number fields.

Over $\R$, we set $\abs{\cdot}_\infty$
to be the usual absolute value (such that $\abs{2}_\infty=2$!).
If $p$ is a prime number, the absolute value over $\Q_p$
is normalized by $\abs{p}_p=1/p$.
These absolute values extend uniquely to any algebraic
extension of $\R$ or $\Q_p$. 

If $F$ is a number field, we denote by $\Val(F)$
the set of places (equivalence classes of valuations) of $F$.
If $v$ is a place of $F$, we will
denote by $\abs{\cdot}_v$ the unique extension of $\abs{\cdot}_\infty$
or $\abs{\cdot}_p$ to $F_v$ (according to $v$ being archimedian or not).
We also set $m_v=e_vf_v$, the product of the ramification
index by the local degree at $v$.
Now, for any $x\in F$ and any valuation $v$ of $\Q$,
\[ \prod_{w|v} \abs{x}_w^{m_w} =\abs{N_{F/\Q}(x)}_v. \]
With these normalizations, we have the \emph{product formula:}
for any $x\in F^*$, 
\[ \prod_{v\in\Val(F)} \abs{x}_v^{m_v}
        = \prod_{v\in\Val(\Q)} \abs{N_{F/\Q}(x)}_v= 1 . \]

Let $X$ be a projective variety over $F$ and $\mathscr L$
a metrized line bundle on $X$. For any $x\in X(F)$,
the (exponential, absolute) \emph{height} of $x$ with
respect to the metrized line bundle $\mathscr L$
is defined by
\[ H_{\mathscr L}(x) = \prod_{v\in\Val(F)} \norm{\mathsf s}_v ^{m_v} (x)
\]
where $\mathsf s $ is any $F$-rational local section of $\mathscr L$,
defined and non-zero at $x$.
The product formula implies that the height doesn't depend on the
choice of $\mathsf s$.

\subsection{Heights on blow-ups}
We keep the notations of the preceeding section.
Moreover we identify $\G_a^n$ with its isomorphic inverse image in $X$
under the blow-up $\pi:X\ra \P^n$.

The metrizations above allow us to define \emph{height functions}
corresponding to complexified divisors $D(\mathbf s)=s_0D_0+\dots+s_rD_r$.
Namely, if $v$ is a place of $F$
and $x=(x_1,\dots,x_n)\in\G_a^n(F_v)$, its exponential local
height is defined by
\[ H_{D(\mathbf s),v}(x) = \prod_{k=0}^r 
       \norm{\mathsf s_{D_k}}_v^{-m_v s_k} (1:x_1:\dots:x_n).
\]
The global height of a point 
 $x\in \G_a^n(\A_F)$ is then the product of all local heights.
This gives a pairing
\[ H : \Pic^G(X)_\C \times \G_a^n(\A_F) \longrightarrow \C^* \]
which is multiplicative as a function on $\Pic^G(X)$
and which is invariant under the
action of the maximal compact subgroup of $\G_a^n(\A_F)$.
Such a pairing had already appeared in the context
of toric varieties.

The invariance of the heights
is a crucial technical ingredient in the proofs
of analytic properties of the height zeta
functions for toric varieties and  
for equivariant compactifications
of $\G_a^n$ considered in  the present paper. 

The ``height zeta function'' is the series
\[ Z(s_0,\dots,s_r)
     = \sum_{x\in\G_a^n(F)} H_{D(\mathbf s)}(x)^{-1} .
\]
Its convergence in some non-empty open subset of $\C^{r+1}$
is a consequence of the
following (well known) lemma.

\begin{lem}\label{prop:zeta-abs}
Let $V$ be a projective variety over a number field $F$
and $({\mathscr L_i})_{1\leq i\leq d}$ a finite number of
ample metrized line bundles on $V$.
For $x\in V(F)$,
define $H(\mathbf s;x)=\prod_{i=1}^d H_{\mathscr L_i}(x)^{s_i}$.
Then  there exists an open non-empty subset $\Omega$ of $\R^d$
such that the series 
\[
Z({\mathbf s})=\sum_{x\in X(F)}H({\mathbf s};x)^{-1}
\]
converges absolutely and uniformly for all ${\mathbf s} \in \C^d$ 
with $\Re({\mathbf s})$ contained in $\Omega$.

Moreover, any other metrization on the $\mathscr L_i$'s
gives the same domain of convergence.
\end{lem}
\begin{proof}
The usual proof of Northcott's theorem establishes
a polynomial bound for the number of rational points
of bounded exponential height.
Hence, the height zeta function of $(\P^n,\mathscr O(1))$
converges for $s\gg 0$.
(There is no need to invoke Schanuel's theorem~\cite{schanuel79}
which gives the precise \emph{asymptotics}.)

Therefore, there are real numbers $\alpha_i$ such that $Z(0,\dots,s_i,\dots,0)$
converges for $\Re(s_i)\geq \alpha_i$.
Now, $Z(\mathbf s)$ converges for any $\mathbf s=(s_1,\dots,s_d)\in \C^{d}$
such that for each $i$, $\Re(s_i)\geq \alpha_i$.
\end{proof}

\subsection{Harmonic analysis on the additive group}
We recal basic facts concerning harmonic
analysis on the group of adelic points $\G_a^n(\A_F)$ 
(cf., for example, \cite{tate67b}). 
For any prime number $p$,
we can view $\Q_p/\Z_p$ as the $p$-Sylow subgroup 
of $\Q/\Z$. This allows us to define a local character $\psi_p$
of $\G_a(\Q_p)$ 
by setting 
\[
\psi_p \,: \, x_p \mapsto \exp(2 \pi i x_p).
\]
At the infinite place of $\Q$ we put 
\[
\psi_{\infty}\colon  x_{\infty} \mapsto \exp(-2 \pi i x_{\infty}), 
\]
(here $x_{\infty}$ is viewed as an element in $\R/\Z$). 
The product of local characters gives a character
$\psi$ of $\G_a(\A_{\Q})$ and, by composition with the trace, 
a character of $\G_a(\A_F)$. For any ${\mathbf a}\in \G_a^n(\A_F)$
we obtain a character $\psi_{\mathbf a} $ of $\G_a^n(\A_F)$
by 
\[
{\mathbf x}  \mapsto  \psi \circ {\rm tr}_{F/\Q}(
\langle \mathbf a , {\mathbf x} \rangle ).
\]
The choice of $\psi$ defines 
a self-duality of $\G_a^n(\A_F)$ (Pontryagin duality).
For $v\in \Val(F)$,
we denote by $\mu_v$ 
the standard normalized local Haar measures
on $\G_a^n(F_v)$ and by $\mu=\prod_v \mu_v$ the self-dual measure 
on $\G_a^n(\A_F)$. The precise normalization can be found in 
(cf. \cite{tate67b} or \cite{manin-p95}, p. 280);
for $F=\Q$, we have $\mu_p(\Z_p)=1$ and $\mu_\infty([0;1])=1$.

For a function $H$ on $\G_a^n(\A_F)$ we denote by $\hat{H}$
its Fourier-transform (with respect to the Haar measure $\mu$)
\[
\hat H:\G_a^n(\A_F)\ra\C,
\qquad \psi \mapsto \int_{\G_a^n(\A_F)}H({\mathbf x})\psi({\mathbf x})
\,\mathrm d\mu(\mathbf x),
\]
whenever the integral converges. 
We shall also use the notation
$\mathrm d\mathbf x$ for $\mathrm d\mu(\mathbf x)$.

We will use the Poisson formula in following form
(cf.~\cite{manin-p95}, p.~280).
\begin{thm}
\label{theo:poisson}
Let $H$ be a continuous 
function on $\G_a^n(\A_F)$ such that
both $H$ and $\hat{H}$ are integrable
and such that 
\[
\sum_{{\mathbf a}\in \G_a^n(F)} H({\mathbf x}+{\mathbf a})
\]
is absolutely and uniformly convergent on compact subsets 
in $\G_a^n(\A_F)/\G_a^n(F)$. Then
\[
\sum_{{\mathbf x}\in \G_a^n(F)} H({\mathbf x})=\sum_{{\mathbf a}\in \G_a^n(F)}
\hat{H}(\psi_{\mathbf a}).
\]
\end{thm}

For $\mathbf s\in\C^{r+1}$ and $\psi\in\G_a^n(\A_F)$,
we shall denote by $\hat H(\mathbf s;\psi)$
the Fourier transform of the height function $H(\mathbf  s;\cdot)^{-1}$
on $\G_a^n(\A_F)$ at the character $\psi$.
It is the product of the local Fourier transforms of
the functions $H_v(\mathbf s;\cdot)^{-1}$ for all $v\in \Val(F)$.

\begin{prop}\label{prop:vanishing}
With the above notations, for all characters
$\psi$ which are non-trivial on the
maximal compact subgroup of $\G_a^n(\A_F)$, we have
\[ \hat{H}({\mathbf s},\psi)=0.\]
\end{prop}
\begin{proof}
This follows from the invariance of the height under the 
action of maximal compact subgroups, see Prop.~\ref{prop:invariance}.
\end{proof}

Consequently, we have a formal identity for the
height zeta function:
\begin{equation} \label{eq:formal}
Z({\mathbf s})=\sum_{\mathbf a\in\G_a^n(\mathfrak o_F)}
    \hat{H}({\mathbf s};\psi_{\mathbf a}).
\end{equation}

The following lemma verifies the two hypotheses
of the Poisson formula~\ref{theo:poisson} concerning $H$.
\begin{lem}
There exists a real $\alpha$ such that for any
$\mathbf s\in\C^{r+1}$ satisfying $\Re(s_0-s_k)\geq\alpha$
and $\Re(s_k)\geq \alpha$, 
and for any compact subset $K$ of $\G_a^n(\A_F)/\G_a^n(F)$,
the series
\[ \sum_{\mathbf a\in\G_a^n(F)} H(\mathbf s;x+\mathbf a) \]
converges absolutely and uniformly for $\mathbf x\in K$.
\end{lem}
\begin{proof}
If $\mathbf s\in\Z^ {r+1}$, the line bundle $D(\mathbf s)$ is ample iff
all $s_k>0$ and $s_0>s_1+\dots+s_r$.
Moreover, the ample line bundles 
$ D= (r+1) D_0+ D_1+\dots+D_r , D+D_1,\dots,D+D_r$
provide a basis of $\Pic(X)_\R$.
Hence, Lemma~\ref{prop:zeta-abs} implies the existence of
a real $\alpha>0$ such that
the series converges absolutely when $\mathbf x=0$,
uniformly for all $\mathbf s\in\C^{r+1}$ such that
$\Re(s_0)>\alpha$, $\Re(s_0-s_k)>\alpha$.

For any $\mathbf x$, the function $H(\mathbf s; \mathbf x+\cdot)$ is another
height function for $D(\mathbf s)$, called ``twisted height'' in our
paper~\cite{chambert-loir-t99}, \S\,2.4, esp.~proposition 2.4.3.
As the convergence of the height zeta function
doesn't depend on metrizatins, this implies the convergence
for any $\mathbf x$. The uniformity for $\mathbf x\in K$
follows from the fact that the height functions can be 
mutually uniformly bounded.
\end{proof}

Now, for the proof of the meromorphic continuation of the height
zeta function it will be sufficient to prove that the
$\hat H$-series on the right-hand side of Equation~\ref{eq:formal},
1\textsuperscript{o}) converges for some
$\Re(s_0)>\alpha$, $\Re(s_0-s_k)>\alpha$  big enough,
and 2\textsuperscript{o}) continues meromorphically.

\subsection{Integrability of local height functions}
The aim of this section is to prove a general result 
concerning the integrability
of local height functions against a measure with singularities.

\begin{prop} \label{prop:integrab}
Let $X$ be a proper smooth variety of dimension $d$ over 
a field $F_v$ which is a finite extension of $\R$ or $\Q_p$.
Fix a finite number of metrized line bundles $\mathscr L_\alpha$ on $X$
together with sections $\mathsf s_\alpha$. Assume that their
divisors $\div(\mathsf s_\alpha)$ are smooth and that their
sum is a divisor $D$ with normal crossings and let $U=X\setminus D$.
Finally, let $\omega\in\Gamma(U,\Omega^d_{X/F_v})$ be a meromorphic
differential form of top degree. We assume that there are
integers $\lambda_\alpha$ such that the divisor of $\omega$
equals
$\sum_\alpha \lambda_\alpha \div(s_\alpha)$.
Denote by $\mathrm d\omega$ the associated measure on $U(F_v)$.

Then, the integral
\[ \int_{U(F_v)} \prod_\alpha \norm{\mathsf s_\alpha}^{m_v r_\alpha} (x)
\, \mathrm d \omega \]
converges if and only if for all $\alpha$, $r_\alpha>\lambda_\alpha-1$.
\end{prop}
\begin{proof}
Using a partition of unity on $X$ for the $F_v$-topology,
we may assume that 
$X$ is a relatively compact open subset $\Omega\subset F_v^d$, with
local coordinates $x_1,\dots,x_d$ and that
the divisor $\sum_\alpha \div(s_\alpha)$ is given by
the equation $x_1\dots x_a=0$
for some integer $0\leq a\leq d$.
The integral is then
\[ I_\Omega = \int_{\Omega}
 \prod_{i=1}^a \abs{x_i}_v^{m_v (r_{\alpha(i)}-\lambda_{\alpha(i)})}
 \exp\big(\sum_{\alpha} h_\alpha(x)\big) \, \mathrm dx_1\dots \mathrm dx_d \]
for some functions $h_{\alpha}$
giving the metrics in our local trivialization
and which are therefore continuous and bounded.

Remark that the integral of $\abs{x}_v^{m_v s}$ over the unit ball of $F_v$
converges if and only if $s>-1$.
The Fubini theorem shows that the integral $I_\Omega$ converges if and only if 
for each $i\in\{1,\dots,a\}$, $r_{\alpha(i)}-\lambda_{\alpha(i)}>-1$.
As any $\alpha$ appears in some chart, the proposition is proved.
\end{proof}

\Subsection{The local Fourier transform in the archimedian case}

When $F_v=\R$ or $\C$,
we want to show that the local Fourier transform of the height
function as a function of $\psi_{\mathbf a}$
decreases rapidly when the norm of $\mathbf a\in F_v^n$ grows to infinity.
The proof proceeds by integration by parts, which requires some estimates.

\begin{lem}
Let $X$ be a smooth projective variety over $F_v$ and
$Z$ be a smooth closed subscheme of $X$.
Let $\partial$ be
a global section of $(\Omega^1_X)^\vee\otimes\mathscr I_Z$,
i.e. a derivation on $X$ vanishing on $Z$.
Denote by $\pi:Y\ra X$ the blow-up of $\mathscr I_Z$.

1) Then the derivation $\partial|_{\pi^{-1}(X\setminus Z)}$ extends uniquely
to a derivation on $Y$.

2) Let $\mathscr E$ be a vector bundle on $X$ equipped
with  a smooth hermitian metric
and $\mathsf s$ a global section of $\mathscr E$ whose
divisor is $Z$. Then the function $\partial \log\norm{\mathsf s}$
extends uniquely to a smooth function on $Y$.
\end{lem}

\begin{proof}
Choose local analytic coordinates on $X$
such that $Z$ is defined by $x_1=\dots=x_a=0$.
Then, $Y$ may be embedded in $\P^{a-1}\times \A^d$ with
coordinates $((t_1:\dots:t_a),(x_1,\dots,x_d))$
and is given there by the equations $t_ix_j=t_j x_i$ for $i\in\{1,\dots,a\}$.
We consider the chart $t_a\neq 0$. Then, local coordinates
on $Y$ are $t_1,t_2,\dots,t_{a-1},x_a,x_{a+1},\dots,x_d$ and
$\pi:Y\ra X$ is given by $x_i=t_i x_a$ if $i<a$.

On $X$, the derivation $\partial$ has the form
\[ \partial = \sum_{i=1}^d h_i \frac{\partial}{\partial x_i}, \]
for some functions $h_i\in (x_1,\dots,x_a)$.
Now, we have to verify that if $i<a$, $\partial t_i$
is a regular function on $Y$. But
\[
\partial t_i  = \partial (x_a/x_i) 
= h_a (x) \frac{1}{x_i} - h_i(x) \frac{x_a}{x_i^2}
= \frac{1}{x_a} h_a (x) t_i - \frac{1}{x_a}  h_i(x) t_i
\in \mathscr O_Y
\]
since 
\[ h_j(x) \in (t_1x_a,\dots,t_{a-1}x_a,x_a)=(x_a). \]

For the statement concerning norms, 
we may fix the coordinates so that 
$ \norm{\mathsf s}^2(x) = \sum_{i=1}^a \abs{x_i}^2 $.
Then, 
\begin{align*}
 \partial\log\norm{\mathsf s}^2 &
= \frac{1}{\abs{x_1}^2+\dots+\abs{x_a}^2}
\big( \sum_{i=1}^a 2x_i h_i(x) \big) \\
& = \frac{1}{\abs{t_1}^2+\dots+\abs{t_{a-1}}^2+1} 
\big(  \sum_{i=1}^a 2 t_i \frac{h_i(x)}{x_a}\big) 
\end{align*}
is regular on $Y$.
\end{proof}

\begin{prop} \label{prop:infinity}
For any compact subset $K\subset\R^{r+1}$ where $H_v(\mathbf s;\cdot)^{-1}$
is integrable, and for any integer $d\geq 1$, there exists a constant
$c(d,K)$ such that for any $\mathbf a\in\C^{n}$ and
any $\mathbf s\in \C^{r+1}$ with $\Re(\mathbf s)\in K$,
\[ \abs{ \hat H_v(\mathbf s;\psi_{\mathbf a})}
        \leq c(d,K)\left( \frac{1+\norm{\Im(\mathbf s)}}{1+\norm{\mathbf a}}
                \right)^d. \]
\end{prop}
\begin{proof}
The 2 preceeding lemmas imply that for any multiindex $\alpha\in\N^n$,
the derivative $\frac{\partial^\alpha}{\partial x^\alpha}
(\log \norm{\mathsf s_D})(x)$ is bounded on $\G_a^n(F_v)$.
Moreover, $\norm{\mathsf s_D}$
tends to $0$ at infinity. We thus may integrate by parts $d$-times.
\end{proof}

\section{Projective space}
\label{sect:pn}

From now on, we work over the field of rational numbers 
$\Q$. It will be clear from the proofs that the case of general
number fields is indeed similar.

This section is included to illustrate our approach
in the simplest example:
we give yet another proof 
of asymptotics for the number of rational points
of bounded height on the standard projective 
space $\P^n$ over the field of rational numbers with the standard
metrization of the line bundle ${\mathcal O}(1)$ given 
by the model $\P^n_{\Z}$ at the finite places and 
by the $L^2$-norms at the archimedian places. 

To keep this section
as self-contained as possible, we reprove
the estimates needed without referring to the
general estimates of the preceeding section.

We will denote by $\A$ the ring of adeles $\A_{\Q}$, by
$p$ a prime number.
We have the normalized valuations
$|\cdot |_{p}$ with $|p|_{p} = p^{-1}$ 
and the usual absolute value $|\cdot |_\infty$.
If $\mathbf a\in\G_a^n(\A)$,
we denote by $\psi_{\mathbf a}$ the corresponding character 
via the identification of $\G_a^n(\A)$ with its Pontryagin dual.

We are interested in the height zeta function
\begin{equation}\label{eqn:zeta-pn}
Z(s)=\sum_{{\mathbf x}\in \G_a^n(\Q)} H({\mathbf x})^{-s}
\end{equation}
where 
$H({\mathbf x}) = H_\infty (\mathbf x) \, \prod_{p} H_{p}({\mathbf x}) $
with
\[
H_{v}({\mathbf x}) := \norm{{\mathbf x}}_v =
 \begin{cases}
(1+\sum_{j=1}^n  \abs{x_j}_v^2)^{1/2} & \text{if $v|\infty$}  \\
\max (1,\max_j |x_j|_v) & \text{if $v$ is finite.} 
\end{cases} 
\]
The series (\ref{eqn:zeta-pn}) 
converges absolutely and uniformly
to a holomorphic function for $\Re(s)\gg 0$.
For all $s$ such that the both sides converge,
we have the Poisson-formula identity (cf. \ref{theo:poisson})
\begin{equation}\label{eqn:poisson-pn}
Z(s)=\sum_{\psi_{\mathbf a}} \hat{H}(s;\psi_{\mathbf a}),
\end{equation}
absolutely. This identity is the starting point
for a meromorphic continuation of $Z(s)$. We now
compute (resp. estimate) the local Fourier transforms.

\begin{lem}\label{lem:pn-0} 
Let $p$ be a prime number.
For all $s$ with $\Re(s)>n$, $H_p(s;\cdot)$ is integrable
on $\Q_p^n$ and its Fourier transform at the trivial character $\psi_0$
is given by
\begin{equation}\label{eqn:pn-0}
\hat{H}_{p}(s;\psi_0)=\frac{1-p^{-s}}{1-p^{-(s-n)}}.
\end{equation}
\end{lem}
\begin{proof}
We decompose the domain of integration $\Q_p^n$
into subdomains 
\[
U(\alpha )
=\big\{
 {\mathbf x}= (x_1,...,x_n) \,;\, \norm{\mathbf x}_p=p^{\alpha }
 \big\},
\]
for $\alpha \ge 1$ and 
\[
U(0)=\big\{
  {\mathbf x}= (x_1,...,x_n) \,;\, \norm{\mathbf x}_p\le 1
     \big\} 
\]
Then 
\begin{align*}
\hat{H}_p(s;\psi_0) & =\int_{U(0)}  
H({\mathbf x})^{-s}\dx +\sum_{\alpha \ge 1} \int_{U({\alpha })} H({\mathbf x})^{-s} \dx,
\\
& =1+\sum_{\alpha \ge 1} p^{-\alpha  s}\cdot\vol(U({\alpha })).
\end{align*}
One has $\vol U(0)=1$ and for $\alpha\geq 1$,
\[ \vol(U(\alpha ))= p^{\alpha  n}\vol (\Z_p^n\setminus (p\Z_p)^n)
 = p^{\alpha n} (1-p^{-n}).  \]
For all $s$ with $\Re(s)>n$, the geometric series
converges absolutely and we obtain
\begin{align*}
\hat{H}_p(s;\psi_0)
	& =1+(1-\frac{1}{p^n})\sum_{\alpha \ge 1}p^{-\alpha  (s-n)}, \\
	& =1 + (1-\frac{1}{p^n})\cdot \frac{1}{p^{s-n}}
		\cdot \frac{1}{1-p^{-(s-n)}}.
\end{align*}
Simplifying, we obtain~\eqref{eqn:pn-0}.
\end{proof}

For all ${\mathbf a}=(a_1,...,a_n)\in \Z^n$ let $S({\mathbf a})$
be the set of all primes which divide all $a_j$.

\begin{lem}\label{lem:pn-q} 
For all ${\mathbf a}\in \Z^n\setminus\{0\}$,
all $s$ with $\Re(s)>n$ and all $p\not\in S({\mathbf a})$
we have
\begin{equation}\label{eqn:q-a}
\hat{H}_p(s;\psi_{\mathbf a})= 1-p^{-s}.
\end{equation}
\end{lem}
\begin{proof}
As above,
we have
\[
\hat{H}_p(s,\psi_{\mathbf a})=1+ 
\sum_{\alpha \ge 1} 
p^{-\alpha  s}\int_{U(\alpha )}\psi_{\mathbf a}({\mathbf x})\dx.
\]
We first integrate over
the set $V(\alpha)$ of $\mathbf x\in\Q_p^n$ such that
$\norm{\mathbf x}\leq p^\alpha$.
\[
 \int_{V(\alpha)} \psi_{\mathbf a}(\mathbf x)\dx
  = p^{\alpha n} \int_{\Z_p^n} \psi_{\mathbf a/p^\alpha}(\mathbf x)\dx.
\]
If $\alpha\geq 1$, as $p$ doesn't divide all the $a_j$,
this is the integral of a non-trivial character on a compact group, hence
$0$. For $\alpha=0$, we get $1$.
Therefore, 
as $V(0)=U(0)$ and 
$U(\alpha)=V(\alpha)\setminus V(\alpha-1)$ for $\alpha\geq 1$,
\[ \int_{U(\alpha)} \psi_{\mathbf a}(\mathbf x)\dx
 = \begin{cases} 0 &\text {for $\alpha\geq 2$} \\
                -1 & \text{for $\alpha=1$.}
    \end{cases}\]
This implies the lemma.
\end{proof}

\begin{lem}\label{lem:estimate-sa-pn}
For all $\epsilon >0$ there exist constants
$c$ and $\delta>0$ such that for all $s$
with $\Re(s)>n+\epsilon $ and all ${\mathbf a}\in \Z^n\setminus\{0\}$
we have the estimate
\begin{equation}\label{eqn:est-sa}
\abs{\prod_{p\in S({\mathbf a})}
\hat{H}_p(s;\psi_{\mathbf a})} \le c\cdot (1+\norm{\mathbf a})^{\delta }.
\end{equation}
\end{lem}

\begin{proof}
In the integral, we replace $\psi_{\mathbf a}$ by 1,
$s$ by $\Re(s)$ and 
use the computation in (\ref{lem:pn-0}). 
For $\Re(s)\geq n+\eps$, we obtain 
\[ \abs{\hat H_p(s;\psi_{\mathbf a})}\leq \frac{1}{1-p^{-\eps}}. \]
If $a$ is a positive integer, we have an inequality
\[ \prod_{p|a} \frac{1}{p^\eps} \ll \ln (1+a) \]
which can be deduced e.g.\ from the Prime Number Theorem.
This gives us equation~\eqref{eqn:est-sa}.
\end{proof}

We now turn to the estimations of the local Fourier
transform for the place at infinity.
For the trivial character, we can---as we could in
the non-archimedian case---explicitely compute the relevant integral:
\begin{lem} \label{lem:pn-infty-0}
For all $s$ with $\Re(s)>n$, $H_\infty(s;\cdot)$
is integrable on $\R^n$ and its Fourier transform at
the trivial character $\psi_0$ is given by
\[ \hat H_\infty(s;\psi_0) = \pi^{n/2} \frac{\Gamma((s-n)/2)}{\Gamma(s/2)}.
\]
\end{lem}

\begin{lem}\label{lem:pn-infty}
For all $\delta >0$ and  all compacts ${K}$ in 
the domain $\Re(s)>n$ there exists a constant 
$c(\delta ,{K})$ such that for all ${\mathbf a}\in \Z^n$ and all 
$s\in {K}$  we have
\[
\abs{\hat{H}_{\infty}(s;\psi_{\mathbf a})}
\le c(\delta ,{K})(1+\abs{\Im s})^\delta (1+\norm{\mathbf a})^{-\delta }
\]
\end{lem}
\begin{proof}
By a unitary change of variables, we may assume $\mathbf a=(\norm
{\mathbf a},0,\dots,0)$. Thus, 
\begin{align*}
\hat{H}_{\infty}(s;\psi_{\mathbf a})
&= \int_{\R^n} (1+\norm x^2)^{-s/2} \exp (-2\pi i \norm {\mathbf a} x_1)\, \dx \\
&= \int_\R \int_{\R^{n-1}}
        (1+\abs{ x_1}^2+\norm{\mathbf x'})^{-s/2}
\exp (-2\pi i \norm a x_1)\, \mathrm dx_1 \, \dx' \\
&= \int_\R (1+\abs{x_1}^2)^{-(s-n+1)/2} \exp(-2\pi i \norm a x_1)\, dx_1
    \int_{\R^{n-1}} \frac{\mathrm d\mathbf y}{(1+\norm {\mathbf y}^2)^{s/2}}
\end{align*}
For any $k>0$, the $k$th derivative of $t\mapsto (1+t^2)^{-s}$ is of
the form $P_k(t)(1+t^2)^{-s-k}$ with $P_k$ a polynomial of degree $k$
whose coefficients are polynomials in $s$.
Thus we can integrate by parts and get for any $k$ an expression
\[ \int_\R (1+t^2)^{-(s-n+1)/2} \exp(-2\pi i \norm a t)\, dt
   = \frac{1}{(\pi i\norm a)^k}
    \int_\R \frac{P_k(t)}{(1+t^2)^{k-\frac{s-n+1}2}}\, dt \]
which imply the lemma.
\end{proof}

\begin{rem}
It follows from the arguments above
that the Fourier transform has polynomial growth in vertical strips.
\end{rem}

\begin{thm}
The series
\[
Z(s)=\sum_{\psi_{\mathbf a}} \hat{H}(s;\psi_{\mathbf a})
\]
converges absolutely and uniformly to a holomorphic
function for $s$ with $\Re(s)>n+1$.
The function $Z(s)$ admits
a meromorphic continuation to the domain $\Re(s)>n $
with exactly one simple pole at $s=n+1$.
The residue at this pole equals
\[
\lim_{s\ra n+1} (s-(n+1))\hat{H}(s;\psi_0)=
\lim_{s\ra n+1} (s-(n+1))\int_{\G_a^n(\A_\Q)}H(s;{\mathbf x})\dx.
\]
\end{thm}

\begin{proof}
Choose a real number $\delta>n$.
From the lemmas above, it follows that there exists
$\delta>0$ such that for any compact $K\in\left]n;+\infty\right[$,
any $\mathbf a\in\Z^n\setminus\{0\}$,
and any real $\delta'>0$,
the product of the local Fourier transforms at the character
$\psi_{\mathbf a}$ converges to a holomorphic function of $s$
which satisfies the inequality
\[ \abs{\hat H(s;\psi_{\mathbf a})}
    \leq c(K) (1+\abs{\Im(s)})^ {\delta+\delta'} (1+\norm{a})^{-\delta'},
\qquad \Re(s)\in K. \]
Hence, the sum over all non-trivial $\psi$ converges absolutely
and locally uniformly to a holomorphic function in the domain $\Re(s)>n$.

At the trivial character, we have, if $\Re(s)>n+1$,
\[ \hat H(s;\psi_0) =
\frac{\zeta(s-n)\Gamma((s-n)/2)}{\zeta(s)\Gamma(s/2)}. \]
This has a simple pole at $s=n$ and extends meromorphically
to the domain $\Re(s)>n$, with no other pole there.
\end{proof}

The identification of the residue and Peyre's Tamagawa constant
in~\cite{peyre95}
is straightforward, granted the meromorphic continuation
of $\hat H(s;\psi_0)$.

\section{Blow-ups of $\P^2$}
\label{sect:p2}

\subsection{Preliminaries}

We continue to work over $\Q$ and we keep the notations
of previous sections. 

Let us consider the projective plane $\P^2$ with 
coordinates $(x_0,x_1,x_2)$ and its Zariski open subset
$U\subset \P^2$ given by $x_0\neq 0$. 
Denote by $X$ the blow-up of $\P^2$ in $r$ distinct points 
$Z_1,\dots,Z_r$ contained in the line at infinity
$Z_0\subset \P^2$ which is given by $x=0$.

For all $k\in\{1,\dots,r\}$, there is a linear form
$\ell_k\in \Z[x_1,x_2]$ with coprime coefficients
such that $Z_k=\mathscr V(x_0,\ell_k)$.
For $k=1,\dots,r$, we denote  by $D_k$ the inverse image 
of $Z_k$ in $X$ and by $D_0$ the strict
transform of the line $Z_0$. 
The variety $X$ is smooth; the anticanonical class is given by
\[
[K_X^{-1}]=3[D_0]+2\sum_{k=1}^{r-1}[D_k].
\]

In the sequel, we shall identify a point $\mathbf x\in\G_a^2$
with the point with homogeneous
coordinates $(1:\mathbf x)$ in $\P^2$ or with its pre-image in
the blow-up.
It follows from the general theory of height functions on blow-ups
that for all $k\in\{1,\dots,r\}$, a local height function for the divisor
$D_k$ at such a point $\mathbf x$ is given by
\[ H_{k,p}(\mathbf x)
=  \frac{\max(1,\norm{\mathbf x}_p)}{\max(1,\abs{\ell_k(\mathbf x)}_p)}
\]
at a finite place $p$, and by an analogous formula where
$\max(1,\cdot)$ is replaced by $\sqrt{1+\cdot^2}$ at the infinite
place.
For $D_0$, we have
\[ H_{0,p} (\mathbf x)= \max(1,\norm{\mathbf x}_p) \prod_{k=1}^r
H_{k,p}^{-1}(\mathbf x) \]
(with the same convention if $v=\infty$).
The global height is given by  
\[ H_k (\mathbf x)= H_{k,\infty} (\mathbf x)\cdot \prod_p H_{k,p} (\mathbf x)\]
and for $\mathbf s=(s_0,\dots,s_r)\in\C^r$, we define
\[ 
H({\mathbf s}; {\mathbf x}):=\prod_{k=0}^{r} H_{k}({\mathbf x})^{s_k}
\]
the global height corresponding to the complexified line bundle
$D(\mathbf s)$.


From~\ref{theo:poisson}, we see that the height zeta function for $X$
has the following formal ``Fourier expansion'':
\[
\sum_{{\mathbf x}\in \Q^2}
H({\mathbf s};{\mathbf x})^{-1} =
\sum_{{\mathbf a}\in \Z^2} \hat{H}({\mathbf s};\psi_{{\mathbf a}})
\]
We have the decomposition
\[
\hat{H}({\mathbf s};\psi_{{\mathbf a}})= 
\hat{H}_{\infty}({\mathbf s};\psi_{{\mathbf a},\infty})\cdot
\prod_p\hat{H}_p({\mathbf s};\psi_{{\mathbf a},p}).
\]
As in the case of $\P^n$, we compute the local Fourier
transforms for almost all places and
estimate them at the remaining bad places.

Let $S$ be the set of primes
of bad reduction of the schematic closure of $\bigcup_k Z_k$
in $\P^2_{\Z}$.
A prime $p$ belongs to $S$ if there exist two linear forms
$\ell_k$ and $\ell_j$ such that
$p$ divides $\det(\ell_k,\ell_j)$.

\subsection{Decomposition of the domain}
Fix a prime $p\not\in S$. We may omit $p$ from the notations
for norms, etc.
Define subsets of $\Q_p^2$ as follows:
\begin{itemize}
\item $U(0)=\Z_p^2$;
\item if $1\leq \beta < \alpha$ and $k\in\{1,\dots,r\}$,
        $U_k(\alpha,\beta)$
	is the set of $\mathbf x\in \Q_p^2$
	such that $\norm{\mathbf x}=p^\alpha$
	and $\abs{\ell_k(\mathbf x)}=p^{\alpha-\beta}$;
\item  if $\alpha\geq 1$ and $k\in\{1,\dots,r\}$,
        $U_k(\alpha)$ is the set of $\mathbf x\in \Q_p^2$ such that
	$\norm{\mathbf x}=p^\alpha$ and $\abs{\ell_k(\mathbf x)}\leq 1$;
\item   if $\alpha\geq 1$, $U(\alpha)$ is the set of $\mathbf x\in\Q_p^2$
	such that $\norm{\mathbf x}=p^\alpha$
	and all $\abs{\ell_j(\mathbf x)}=p^\alpha$.
\end{itemize}

As $p\not\in S$, these sets furnish a partition of $\Q_p^2$.
This decomposition is well adjusted to our local heights
since they are constant on each subset:
\begin{itemize}
\item on $U(0)$, all $H_k$'s are $1$;
\item on $U_k(\alpha,\beta)$, $H_k=p^{\beta}$,
the other $H_j$ with $j\geq 1$ are $1$ and $H_0=p^{\alpha-\beta}$;
\item on $U_k(\alpha)$, $H_k=p^\alpha$, the
other $H_j$ are $1$; 
\item on $U(\alpha)$, $H_0=p^\alpha$ and all other are $1$.
\end{itemize}
In other words,
\begin{equation}
H(\mathbf s;\mathbf x) =
\begin{cases}
  1
	& \text{if $\mathbf x\in U(0)$;} \\
  p^{\alpha s_0}
	& \text{if $\mathbf x\in U(\alpha)$, $\alpha\geq 1$;} \\
  p^{\alpha s_0 + \beta (s_k-s_0)}
	& \text{if $\mathbf x\in U_k(\alpha,\beta)$,
		$1\leq \beta<\alpha$;} \\
  p^{\alpha s_k}
	& \text{if $\mathbf x\in U_k(\alpha)$.}
\end{cases}
\end{equation}

The table on p.~\pageref{table:sumup}
summarizes this information.

\subsection{Some integrals of characters}
We identify $\mathbf a\in\Z^2$ with the linear form
it defines on $\G_a^2$ as well as with the character
$\psi_{\mathbf a}$ of $\G_a^2(\A_\Q)$ it determines.

We will say that a character is \emph{generic} if
$\mathbf a$ is not proportional to any $\ell_k$.
A non-trivial character is special if it is proportional to some
(necessarily unique) $\ell_k$.

If $\mathbf a$ is generic, $S(\mathbf a)$ is the set
of primes $p$ such that $p$ divides some determinant
$\det(\ell_k,\mathbf a)$.

If $\mathbf a$ is special for $\ell_k$, $S(\mathbf a)$
is the set of primes such that $p$ divides some
determinant $\det(\ell_j,\mathbf a)$ for $j\neq k$.

Note that $S\subset S(\mathbf a)$ for any non-trivial $\mathbf a$
and that if $p|\mathbf a$, then $p\in S(\mathbf a)$.

We now compute the integral of $\psi_{\mathbf a}$
over the subsets defined in the previous subsection,
at least in the cases $\mathbf a=0$, $\mathbf a$ special
and $\mathbf a$ generic.

Remark that for any $\mathbf a$,
$\int_{U(0)} \psi_{\mathbf a}(\mathbf x) \dx =1$.

\begin{lem}[Trivial character] Let $p$ be a prime not in $S$.
Then, \addtocounter{equation}{-1}
\begin{subequations}
\begin{align}
\label{vol:uk(a,b)}
\vol U_k(\alpha,\beta) & = p^{2\alpha-\beta} \frac{(p-1)^2}{p^2}; \\
\label{vol:uk(a)}
\vol U_k(\alpha) & = p^\alpha \frac{p-1}p ; \\
\label{vol:u(a)}
\vol U(\alpha) & = p^{2\alpha} \frac{(p-1)(p+1-r)}{p^2}.
\end{align}
\end{subequations}
\end{lem}

\begin{lem}[Generic character]
Let $\mathbf a$ be a generic character and $p\not\in S(\mathbf a)$.
Then, \addtocounter{equation}{-1}
\begin{subequations}
\begin{align}
\label{gen:uk(a,b)}
\int_{U_k(\alpha,\beta)} \psi_{\mathbf a} & =  0 \\
\label{gen:uk(a)}
\int_{U_k(\alpha)} \psi_{\mathbf a} & = 
   \begin{cases} -1 & \text{if $\alpha=1$}\\ 0 &\text{else;} \end{cases}\\
\label{gen:u(a)}
\int_{U(\alpha)} \psi_{\mathbf a} & =
   \begin{cases} -1+r & \text{if $\alpha=1$}\\ 0 &\text{else;} \end{cases}
\end{align}
\end{subequations}
\end{lem}

\begin{lem}[Special character]
Let $\mathbf a$ a character which is special
for $\ell_k$. If $p\not\in S(\mathbf a)$ and $j\neq k$, one has
\addtocounter{equation}{-1}
\begin{subequations}
\begin{align}
\label{spec:uj(a,b)}
\int_{U_j(\alpha,\beta)} \psi_{\mathbf a} & =  0 \\
\label{spec:uj(a)}
\int_{U_j(\alpha)} \psi_{\mathbf a} & = 
   \begin{cases} -1 & \text{if $\alpha=1$}\\ 0 &\text{else;} \end{cases}\\
\label{spec:uk(a,b)}
\int_{U_k(\alpha,\beta)} \psi_{\mathbf a} & = 
   \begin{cases} -p^\alpha\frac{p-1}p & \text{if $\beta=\alpha-1$}\\
		 0 &\text{else;} \end{cases}\\
\label{spec:uk(a)}
\int_{U_k(\alpha)} \psi_{\mathbf a} & =
        p^\alpha \frac{p-1}p\\
\label{spec:u(a)}
\int_{U(\alpha)} \psi_{\mathbf a} & =
   \begin{cases} -(p+1-r) & \text{if $\alpha=1$}\\ 0 &\text{else;} \end{cases}
\end{align}
\end{subequations}
\end{lem}

\begin{proof}
We prove the three lemmas simultaneously.
By a unitary change of variables, we may assume that $\ell_k(\mathbf x)=x_1$.
Then one has
\[ U_k(\alpha,\beta)=p^{\beta-\alpha}\Z_p^* \times p^{-\alpha} \Z_p^* \]
and
\[ U_k(\alpha) = \Z_p \times p^{-\alpha}\Z_p^*, \]
hence their volumes as in formulas~\eqref{vol:uk(a,b)} and~\eqref{vol:uk(a)}.

If $p$ doesn't divide $\det(\ell_k,\mathbf a)$, we may change
variables and even assume that $\mathbf a=(0,1)$.
Then,
\[ \int_{U_k(\alpha,\beta)} \psi_{\mathbf a}(\mathbf x)\dx
 = p^{2\alpha-\beta} \frac{p-1}p \int_{\Z_p^*} \exp(2\pi i u/p^\alpha)\, du
\]
and the last integral has already been calculated when we studied
the case of $\P^n$ : one finds $0$ if $\alpha\geq 2$ and
$-1/p$ if $\alpha=1$. But $\alpha>\beta\geq 1$, so $\alpha\neq 1$.
This proves formulas~\eqref{gen:uk(a,b)} and~\eqref{spec:uj(a,b)}.

Similarly, 
\[ \int_{U_k(\alpha)} \psi_{\mathbf a}(\mathbf x)\dx
= p^{\alpha} \int_{\Z_p^*} \exp(2\pi i u/p^\alpha)\, du
\]
is $-1$ for $\alpha=1$ and $0$ else.
Formulas~\eqref{gen:uk(a)} and~\eqref{spec:uj(a)} are
therefore proved.

We now treat the case of a character $\mathbf a$ which is special
for $\ell_k$.
A unitary change of variables allows to assume $\ell_k(\mathbf x)=x_1$
and $\mathbf a=(1,0)$.
Then,
\[ \int_{U_k(\alpha,\beta)} \psi_{\mathbf a}(\mathbf x)\dx
= p^{2\alpha-\beta}\frac{p-1}p \int_{\Z_p^*} \exp(2\pi i x/p^{\alpha-\beta})
\]
is $0$ if $\alpha-\beta\neq 1$ and is equal to
\[ p^{2\alpha-\alpha+1} \frac{p-1}p \frac{(-1)}p = - p^\alpha \frac{p-1}p \]
if $\alpha=\beta-1$, as stated in~\eqref{spec:uk(a,b)}.
Equation~\eqref{spec:uk(a)} follows from
\[ \int_{U_k(\alpha)} \psi_{\mathbf a}(\mathbf x)\dx
= \int_{\Z_p \times p^{-\alpha}\Z_p^*} \exp(2\pi i x_1)\dx
= p^\alpha\frac{p-1}p.
\]
To compute the volume of $U(\alpha)$, it is useful to remark that
$U(\alpha)$ is $p^{-\alpha}$ times the complementary subset
in $\Z_p^2$ of $1+(p-1)r$ disjoint balls of radius $1/p$.
Therefore,
\[ \vol U(\alpha)=p^{2\alpha} \big (1 -\frac{1+(p-1)r}{p^2}\big) 
        = p^{2\alpha} \frac{(p-1)(p+1-r)}{p^2}, \]
as in formula~\eqref{vol:u(a)}.

If $p\nmid\mathbf a$, remark that the integral of $\psi_{\mathbf a}$
over $p^{-\alpha}(\Z_p^2\setminus p\Z_p^2)$
is $-1$ for $\alpha=1$ and $0$ for $\alpha\geq 2$.
We now need to substract the integrals
over the disjoint subsets  $U_k(\alpha,\beta)$ and $U_k(\alpha)$.

For a generic character, one gets $0$ if $\alpha\geq 2$ and
$-1+r$ if $\alpha=1$; 
this establishes formula~\eqref{gen:u(a)}.
Finally, if $\mathbf a$ is special for $\ell_k$, one has
\[ \int_{U(1)} \psi_{\mathbf a}(\mathbf x)\dx
= -1 + (r-1) - (p-1) = - (p+1-r) \]
and $ \int_{U(\alpha)}\psi_{\mathbf a}=0$ for $\alpha\geq 2$,
as claimed in~\eqref{spec:u(a)}.
\end{proof}

\subsection{The local Fourier transform  at $\psi_0$}
We still assume $p\not\in S$ and compute the local Fourier
transform at the trivial character $\psi_0$.
By the general result~\ref{prop:integrab}
on Fourier transforms of height functions,
$H(\mathbf s;\cdot)^{-1}$ is integrable on $\G_a^n(\Q_p)$
as soon as $\Re(s_0)>2$ and all $\Re(s_k)>1$.
We then have:
\[
\hat H(\mathbf s;\psi_0) = \int_{U(0)} +
\left(\sum_{k=1}^r \sum_{1\leq\beta<\alpha} \int_{U_k(\alpha,\beta)}
+ \sum_{\alpha=1}^\infty \int_{U_k(\alpha)} \right)
+\sum_{\alpha=1}^\infty \int_{U(\alpha)}
\]
and we compute each sum separately.
The integral over $U(0)$ is $1$.
Now, for a fixed $k$, the integral over all $U_k(\alpha,\beta)$ is
\begin{align*}
\sum_{1\leq\beta<\alpha} \int_{U_k(\alpha,\beta)}
&=
    \frac{(p-1)^2}{p^2} \sum_{1\leq\beta<\alpha}
        p^{-\alpha s_0} p^{-\beta(s_k-s_0)} p^{2\alpha} p^{-\beta} \\
&=    \frac{(p-1)^2}{p^2} \sum_{\beta=1}^\infty p^{-\beta(s_k-s_0+1)}
            \sum_{\alpha=\beta+1}^\infty p^{-\alpha(s_0-2)} \\
&= \frac{(p-1)^2}{p^2} \sum_{\beta=1}^\infty p^{-\beta(s_k-s_0+1)}
            p^{-\beta(s_0-2)} \frac{1}{p^{s_0-2}-1}  \\
&=    \frac{(p-1)^2}{p^2} \frac{1}{p^{s_0-2}-1} 
        \sum_{\beta=1}^\infty p^{-\beta(s_k-1)}\\
&= \frac{(p-1)^2}{p^2} \frac{1}{p^{s_0-2}-1}
            \frac{1}{p^{s_k-1}-1}
\end{align*}

The sum over all $U_k(\alpha)$ ($k$ fixed) equals
\[ \sum_{\alpha=1}^\infty \int_{U_k(\alpha)}
=    \frac{p-1}p \sum_{\alpha=1}^\infty p^{-\alpha s_k} p^\alpha 
=    \frac{p-1}p \sum_{\alpha=1}^\infty p^ {-\alpha(s_k-1)} 
=    \frac{p-1}p \frac{1}{p^{s_k-1}-1}
\]
Finally, the sum over all $U(\alpha)$ is
\begin{align*}
\sum_{\alpha=1}^\infty \int_{U(\alpha)}
&= \frac{(p-1)(p+1-r)}{p^2}
	\sum_{\alpha=1}^\infty p^{-\alpha s_0} p^{2\alpha} 
= \frac{(p-1)(p+1-r)}{p^2} \sum_{\alpha=1}^\infty p^{-\alpha(s_0-2)} \\
&= \frac{(p-1)(p+1-r)}{p^2} \frac{1}{p^{s_0-2}-1} 
=  \frac{(p-1)(p+1-r)}{p^2} \frac{1}{p^{s_0-2}-1} 
\end{align*}

Putting all this together,
we have
\begin{align}
\hat H(\mathbf s;\psi_0) 
&= 1+ \frac{p-1}{p^2} \frac{1}{p^{s_0-2}-1}
           \sum_{k=1}^r \frac{1}{p^{s_k-1}-1}
            \left( (p-1) + p (p^{s_0-2}-1)\right) \notag \\
&\qquad
        + \frac{(p-1)(p+1-r)}{p^2} \frac{1}{p^{s_0-2}-1} \notag \\
&= 1+ \frac{p-1}{p^2} \frac{p^{s_0-1}-1}{p^{s_0-2}-1} 
            \sum_{k=1}^r \frac{1}{p^{s_k-1}-1} 
    + \frac{(p-1)(p+1-r)}{p^2} \frac{1}{p^{s_0-2}-1} \notag \\
&= 1+ \frac{p-1}{p^2} \frac{1}{p^{s_0-2}-1}
        \left( (p+1-r)+ \sum_{k=1}^r \frac{p^{s_0-1}-1}{p^{s_k-1}-1} \right)
\notag \\
&= 1+ \frac{p^2-1}{p^{s_0}-p^2}
        + \frac{p-1}{p^{s_0}-p^2}
             \sum_{k=1}^r \frac{p^{s_k-1}-p^{s_0-1}}{p^{s_k-1}-1} 
\label{fourier:1}
\end{align}     

We remark that if $(s_0,s_1,\dots,s_r)=(3,2,\dots,2)$, corresponding
to the anticanonical class $K_X^{-1}$, this yields
\[
\hat H(K_X^{-1},\psi_0) = 1+ \frac{p^2-1}{p^3-p^2} 
	+ r \frac{p-1}{p^3-p^2}  \frac{p^2-p}{p-1}
	= 1+ \frac{r+1}{p} + \frac{1}{p^2}
        = \frac{1}{p^2} \# X(\F_p),
\]
the expected local density at $p$.

\subsection{The local Fourier transform at a generic character}
Let $\mathbf a$ be a generic character and $p\not\in S(\mathbf a)$.
In that case, the summation is easier.
The integrals over $U_k(\alpha,\beta)$ are $0$,
as are the integrals over $U_k(\alpha)$ or $U(\alpha)$ if $\alpha\geq 2$.
Therefore
\begin{equation}
\hat H(\mathbf s;\psi_{\mathbf a})
= 1 -  \sum_{k=1}^r p^{-s_k} + (r-1) p^{-s_0}.
\label{fourier:gen}
\end{equation}
For $K_X^{-s}$, this specializes to
\[ \hat H(K_X^{-s},\psi_{\mathbf a})
     = 1 -r p^{-2s} + (r-1)p^{-3s}. \]

\subsection{The local Fourier transform at a special character}
If $\mathbf a$ is special for $\ell_k$ and $p\not\in S(\mathbf a)$,
it behaves as if it
were generic for the other $\ell_j$. Besides
$U(0)$, $U(1)$ and $U_j(1)$ for $j\neq k$,
remain the integrals
over $U_k(\alpha,\alpha-1)$ for $\alpha\geq 2$ and the one
over $U(\alpha)$ for $\alpha\geq 2$.
\begin{align}
\hat H(\mathbf s;\psi_{\mathbf a}) &=
 1- \sum_{j\neq k} p^{-s_j} + (r-p-1) p^{-s_0}  \notag \\
& \qquad {} + \sum_{\alpha=1}^\infty p^{\alpha}\frac{p-1}p p^{-\alpha s_k}
    - \sum_{\alpha\geq 2} p^{\alpha}\frac{p-1}p p^{-\alpha s_0-(\alpha-1)(s_k-s_0)}  \notag \\
&= 1- \sum_{j\neq k} p^{-s_j} + (r-p-1) p^{-s_0} \notag \\
& \qquad {} +  \frac{p-1}p \frac{1}{p^{s_k-1}-1}
		- \frac{p-1}p p^{s_k-s_0} \frac{p^{-s_k-1}}{p^{s_k-1}-1} \notag \\
&= 1- \sum_{j\neq k} p^{-s_j} + (r-p-1) p^{-s_0} 
 +  (p-1) (1-p^{1-s_0})  \frac{1}{p^{s_k}-p}
\label{fourier:spe-k}
\end{align}
For $K_X^{-s}$, one has
\[  \hat H(K_X^{-s},\psi_{\mathbf a})
        = 1- (r-1)p^{-2s} + (r-p-1)p^{-3s} + (p-1)
			\frac{1-p^{1-3s}}{p^{2s}-1}.
\]

\subsection{Bad reduction, general estimates}
If $p\in S(\mathbf a)$, then
the previous analysis doesn't say anything about
the behaviour of $\hat H(\mathbf s;\psi_{\mathbf a})$.
However, for any compact contained
in the domain of integrability of the
height function, there is a uniform
estimate
\[ \abs{\hat H_p(\mathbf s;\psi_{\mathbf a})} \leq C \]
where the constant $C$ doesn't depend on $p\in S$.
For $p\in S(\mathbf a)\setminus S$,
we replace $\psi_{\mathbf a}$ by $1$ and insert
the estimates for the trivial character.

It follows that 
\[ \prod_{p\in S(\mathbf a)} \abs{\hat H_p (\mathbf s;\psi_{\mathbf a})}
\leq C' (1+\norm{a})^\kappa\]
for some constant $C'$ and some $\kappa>0$.

\subsection{Meromorphic continuation}
We split the sum over all characters in $r+2$ parts:
the trivial character is treated separately; then the generic characters;
then the characters which are special for $\ell_k$, $k$ varying
from $1$ to $r$:
\[ \sum_{\mathbf a} \hat H(\mathbf s;\psi_{\mathbf a})
 = \hat H(\mathbf s;\psi_0)
+ \sum_{\text{$\mathbf a$ generic} } \hat H(\mathbf s;\psi_{\mathbf a})
+ \sum_{k=1}^r \sum_{\text{$\mathbf a$ special for $\ell_k$}}
\hat H(\mathbf s;\psi_{\mathbf a}).\]
Let $Z_0$, $Z_{\text{gen}}$ and $Z_{k}$ (for $1\leq k\leq r$)
be the functions defined by the corresponding series.

Each global Fourier transform at a generic character
defines a holomorphic function of $\mathbf s$ in the domain
$\Re(s_0)>2$ and $\Re(s_k)>1$ for all $1\leq k\leq r$.
Moreover, the estimate at infinity~\ref{prop:infinity}
ensures that the sum over
all generic characters converges locally uniformly.
Therefore, $Z_{\text{gen}}$ extends to a holomorphic
function in that domain.

For the characters which are special for $\ell_k$, the product
of the local Fourier transform defines a holomorphic function 
of $\mathbf s$ in the domain
$\Re(s_0)>2$, $\Re(s_j)>1$ if $j\neq k$ and $\Re(s_k)>2$.
Therefore $Z_k$
extends to a meromorphic function in the domain $\Re(s_0)>2$,
$\Re(s_j)>1$ and $\Re(s_k)>1$ with a simple pole along
the hypersurface $s_k=2$.

And finally,
for the trivial character, we have absolute convergence of
the Euler product for
$\Re(s_0)>3$ and $\Re(s_k)>2$ for $1\leq k\leq r$,
and $Z_0$ has meromorphic continuation to the domain
$\Re(s_0)>2$, $\Re(s_k)>1$,
with principal part $1/(s_0-3)(s_1-2)\dots(s_r-2)$.

The estimates of Prop.~\ref{prop:infinity}
as well as standard estimates for the growth of the Riemann zeta function
in vertical strips imply that (away from poles) $Z(\mathbf s)$ has
polynomial growth in vertical strips.

Therefore, we have proven the following theorem:
\begin{thm}
The height zeta function $Z(\mathbf s)$ converges in
the domain $\Re(s_0)>3$, $\Re(s_k)>2$. Moreover, there
exists a holomorphic function $g$ 
in the domain $\Re(s_0)>2$, $\Re(s_k)>1$
such that
\[ Z(\mathbf s) = g(\mathbf s) \frac{1}{(s_0-3)(s_1-2)\dots(s_r-2)} \]
and $g(K_X^{-1})\neq 0$.
Moreover, $g$ has polynomial growth in vertical strips.
\end{thm}

\begin{cor} \label{coro:last}
There exists a polynomial $P_X$ of degree $r$ such that
for all $\alpha>2/3$, the number of points of $\G_a^2(\Q)\subset X(\Q)$
of anticanonical height $\leq H$ satisfies
\[ N(U,K_X^{-1},B) = B P_X(\log B) + O(B^{\alpha}). \]
Moreover, if $\tau(K_X)$ denotes
the Tamagawa number, the leading coefficient of $P_X$
is equal to 
\[ \frac{1}{r!} \frac{\tau(K_X)}{3\cdot 2^r}, \]
as predicted by Peyre's refinement of Manin's conjecture.
\end{cor}

\begin{rems}
1\textsuperscript{o}) Our theorem implies a similar asymptotic formula
for arbitrary line bundles $\mathscr L$, provided their class belongs
to the interior of the cone of effective divisors.

2\textsuperscript{o}) It is curious to note that we now have asymptotics
for blow-ups of $\P^2$ in any number of points \emph{on a line}, 
and in three points in general position, while
the asymptotics for blow-ups of $\P^2$ in $4$ points 
\emph{in general positition} are still unknown.
\end{rems}

\begin{table}[p]
\doublerulesep\arrayrulewidth
\tabcolsep2\tabcolsep
\centering
\def\arraystretch{1.6}
\begin{tabular}{ccccc}
\hline\hline
 & $U(0)$ &  $ U_k(\alpha,\beta)$ & $U_k(\alpha)$ &  $U(\alpha)$ \\
\hline\hline
volume 
  & $1$ & $p^{2\alpha-\beta}\frac{(p-1)^2}{p^2} $
	 & $p^\alpha\frac{p-1}p $ & $p^{2\alpha}\frac{(p-1)(p+1-r)}{p^2}$ \\
\hline
$H_0$ & $1$ & $p^{\alpha-\beta}$ & $1$ & $p^\alpha$ \\
$H_j$ ($j\neq k$) & $1$ & $1$ & $ 1 $ & 1 \\
$H_k$ & $1$ & $ p^\beta $ & $p^\alpha $ & 1 \\
$H(\mathbf s;\cdot)$ & $1$ & $p^{\alpha s_0+\beta(s_k-s_0)}$
	& $p^{\alpha s_k}$ & $p^{\alpha s_0}$  \\
\hline
\multicolumn{5}{c}{Integrals of a generic character $\psi_{\mathbf a}$}\\
\hline
$\alpha=1$ & $1$ & &  $-1$ & $-1+r$ \\
$\alpha\geq 2$ & $1$ & $0$ & $0$ & $0$ \\
\hline
\multicolumn{5}{c}{%
Integrals of a character $\psi_{\mathbf a}$ special for $\ell_k$}\\
\hline
any $\alpha$ & $1$ &  & $p^\alpha \frac{p-1}p $ &  \\
$\alpha=1$ &  &  &  & $ - (p+1-r)$  \\
$\alpha\geq 2$ & & & $0$ \\
$\alpha=\beta-1$ &  & $-p^\alpha\frac{p-1}p$ &  &  \\
\hline\hline
\end{tabular}
\bigskip
\caption{} 
\label{table:sumup}
\end{table}

\vskip 0pt plus .3\textheight
\penalty-200
\vskip 0pt plus -.3\textheight

\def\noop#1{}
\providecommand{\bysame}{\leavevmode ---\ }
\providecommand{\og}{``}
\providecommand{\fg}{''}
\providecommand{\smfandname}{et}
\providecommand{\smfedsname}{\'eds.}
\providecommand{\smfedname}{\'ed.}
\providecommand{\smfmastersthesisname}{M\'emoire}
\providecommand{\smfphdthesisname}{Th\`ese}

\end{document}